\documentclass[a4paper,twoside]{article}

\usepackage{epsfig}
\usepackage{subcaption}
\usepackage{calc}
\usepackage{amssymb}
\usepackage{amstext}
\usepackage{amsmath}
\usepackage{amsthm}
\usepackage{amsfonts}
\usepackage{multicol}
\usepackage{pslatex} 
\usepackage{apalike}
\usepackage{tikz-network}
\usetikzlibrary{positioning}

\newtheorem{assumption}{Assumption}[section]

\usepackage{SCITEPRESS}     % Please add other packages that you may need BEFORE the SCITEPRESS.sty package.

\begin{document}

\title{Simulation study for the comparison of power flow models for a line distribution network with stochastic load demands}

\author{\authorname{Christianen, M.H.M.\sup{1}\orcidAuthor{0000-0002-9611-500X}, Vlasiou, M.\sup{1,2}\orcidAuthor{0000-0002-0457-2925} and Zwart, B.\sup{1,3}}
\affiliation{\sup{1}Eindhoven University of Technology, Eindhoven, The Netherlands}
\affiliation{\sup{2}University of Twente, Enschede, The Netherlands}
\affiliation{\sup{3}Centrum Wiskunde \& Informatica, Amsterdam, The Netherlands}
\email{\{m.h.m.christianen, m.vlasiou\}@tue.nl, bert.zwart@cwi.nl}
}

\keywords{Electric Vehicle Charging; Power Flow Models; Bandwidth-sharing Networks}%The paper must have at least one keyword. The text must be set to 9-point font size and without the use of bold or italic font style. For more than one keyword, please use a comma as a separator. Keywords must be titlecased.}

\abstract{We use simulation to compare different power flow models in the process of charging electric vehicles (EVs) by considering their random arrivals, their stochastic demand for energy at charging stations, and the characteristics of the electricity distribution network. We assume the distribution network is a line with charging stations located on it. We consider the \emph{Distflow} and the \emph{Linearized Distflow} power flow models and we assume that EVs arrive at the network with an exponential rate, have an exponential charging requirement, and that voltage drops on the distribution network stay under control. We provide extensive numerical results investigating the effect of using different power flow models on the performance of the network.}

\onecolumn \maketitle \normalsize \setcounter{footnote}{0} \vfill

\section{\uppercase{Introduction}}
\label{sec:introduction}
%\begin{enumerate}
%    \item Grid capacity problems: real life.
%    \item Congestion in energy networks: model.
%    \item Nonlinear and linear approximations in the driving resource allocation problem
%    \item Goal: empirically show how these approximations compare
%    \item Overview paper
%\end{enumerate}
In recent years, the growing electricity consumption, the active adoption of renewable generation, and the energy transition result in congestion in the electricity network. On one side, more companies use electricity for their production, more houses are heated with heat pumps and more people drive electric cars. On the other side, companies and citizens are generating more and more electricity from wind and sun, which they mostly feed back to the electricity network. This causes network capacity problems, or in other words, congestion. This is illustrated in \cite{Hoogsteen2017}, where the impact of the energy transition on a real electricity network is evaluated. Here, the authors showed that charging a small number of EVs is enough to cause a blackout in a neighborhood. Therefore, it is imperative to study the performance of the network under different power flow models, since these models are used for the design and control of the network.

To model such congestion, we use stochastic processes in which the dynamics are driven by power allocation problems. The stochasticity comes from the random behavior of citizens and companies in consuming and generating electricity, while the dynamics depend on how electricity is allocated among all users. The allocation of electricity or power depends on protocols by network operators but is constrained by physical laws and network constraints. 

Irrespective of the allocation of power, it is important to respect the physical laws of the network and its constraints. In an electricity network, an important constraint is the requirement of keeping voltage losses, or in other words, the voltage drop, on a cable in the network under control. These voltage losses are caused by the physical properties of the cables in the network. Keeping the voltage losses under control ensures that every user in the network receives safe and reliable power at a voltage that is within some standard range, which varies from one country to another \cite{Kerstinga}. For example, according to Dutch law, the voltage drop in a distribution network, a small part of the electricity network, is not allowed to be more than 4.5\% \cite{VanWestering2020}.

In this paper, we consider the specific process of charging EVs in a neighborhood such that the voltage drop in the distribution network stays under control. In this process, the stochasticity comes from random arrivals (at parking lots with charging stations) and charging requirements of EVs, while the power allocated to each EV in the network depends on the number of EVs and the corresponding location of EVs that are currently charging in the network. We model this process as a queue, with EVs representing \emph{jobs}, charging stations classified as \emph{servers}, and the service being delivered as the power supplied to EVs, constrained by physical laws and network constraints. The particular queuing model that we employ can be seen as falling under a general class of queuing networks called \emph{bandwidth-sharing networks}. To model the physical laws and network constraints, we use two approximations of the alternating current (AC) power flow equations \cite{Molzahn2019}, i.e.; we study the \emph{Distflow} and a linearized version of the Distflow model, by ignoring some power losses, that is called the \emph{Linearized Distflow} \cite{Baran1989,Baran1989a}. 

In this simulation study, we assess the accuracy and effectiveness of the Linearized Distflow model compared to the Distflow model. Our goal is to compare the different power flow models on the performance of the stochastic process of EV charging, by the mean number of EVs and the mean charging time of each EV in the network. Furthermore, we gain insights into the behavior of the network by including variability in the distribution of the arrival rates to different parking lots.

The major insights we gain from analyzing such a network are summarized as follows. We observe that the performance of the Linearized Distflow model is comparable to the Distflow model, i.e. the mean number of EVs and the mean charging time of an EV under both power flow models are similar and the relative difference between critical arrival rates (the specific arrival rates under both power flow models for which the mean number of EVs and the mean charging time of an EV explode) is below 5\%. Thus, our first result provides more evidence that using the much simpler Linearized Distflow model is a valid and accurate approximation of the Distflow model, even if the system is highly heterogeneous. Namely, in our numerical examples, we consider cases where one station has almost all the load of the system. Even for such heterogeneous cases, the performance of the network is the same under the Linearized Distflow and the Distflow model. In other words, we do not lose much accuracy from a performance perspective by ignoring some power losses. The second major conclusion is rather surprising. It is very well known from queuing theory that variability in the network causes worse performance. However, what is surprising in this case is that the network does not perform symmetrically under the same loads. If the load of an individual parking lot is way larger than the other loads of the other parking lots, the performance of the network is different from the performance of the network if the (same) largest load is put on another parking lot.

The structure of the paper is as follows. In Section \ref{sec:literature_review}, we provide a literature review on work that has been done on stability results for EV-charging (from an electrical engineering viewpoint), work that has been done on stability in bandwidth-sharing networks (from a mathematical viewpoint), and a comparison between the Linearized Distflow model and the Distflow power flow model. In Section \ref{sec:problem_description}, we provide a detailed model description. In particular, we introduce the queuing network, the constraints and assumptions of the electrical distribution network, and the power flow models we consider. In Section \ref{sec:numerical_study}, we present several numerical experiments showing the accuracy and effectiveness of the Linearized Distflow model and the effect of including variability in the distribution of the arrival rates to the performance of the network.

\section{\uppercase{Literature Review}}\label{sec:literature_review}
First, we discuss literature related to stability for EV charging and stability in bandwidth-sharing networks. Second, we provide literature that compares the Distflow model and the Linearized Distflow model.

The aforementioned concepts of stability are different from each other. The term stability from an electrical engineering perspective relates, for example to situations where the fluctuations of voltages or power losses in the network are too big, the transformers, lines, and cables are overloaded, or where peak demands are too high. The term stability, from a queuing perspective, means the stability of the queuing model and is defined as the positive recurrence of a Markov process. Informally speaking, stability means the ability of all queues to complete the service of all jobs, without the number of outstanding building up infinitely. In our simulations, we consider a finite number of parking spaces at all parking lots, which implies that the queuing model is always stable. Therefore, instead of the stability of the network, we investigate the performance of the network in terms of the mean number of EVs and the mean charging time of an EV in the network but still discuss the literature on stability results. %After paragraph \ref{par:stab_results_EV}, whenever we write \emph{stability}, we mean the stability of the queuing network, unless stated otherwise. 

\paragraph{Stability results  for EV-charging.}\label{par:stab_results_EV}
The literature on stability results for EV charging is limited to numerical experiments. An early paper on stability analysis in EV-charging is by Huang et al \cite{Huang2013}. The authors present a new quasi-Monte Carlo stability analysis method to assess the dynamic effects of plug-in electric vehicles in power systems. They conclude that improvements in stability control are worth further study since the number of EVs is growing. Other simulation studies are conducted to obtain stability conditions. In \cite{DeHoog2014}, the authors explore the constraint of requiring a minimum voltage to charge EVs throughout a network and demonstrate that the physical locations of individual demand and generation of energy play a significant role in determining whether voltages throughout the network remain within required limits or not. Similarly, in \cite{Zhang2016}, the impact of charging EVs on the voltage stability of the distribution network is simulated and analyzed. The simulation results show that the voltage stability is related to the individual loads, total load in the network, and physical properties of the network. In \cite{Ul-Haq2015}, simulations are performed on another test network. For this network, different charging strategies are implemented and it is shown that for some scenarios this can cause significant voltage instability. Last, in \cite{Deb2018}, the authors perform a numerical study on a specific test system, where they investigate the impact of a single EV charging station on the voltage stability, power losses, and economic losses of the distribution network. Here, it is also observed that the location of the EV charging station is important in the smooth operation of the grid.

\paragraph{Stability conditions for queuing networks.} 
The literature of the stability of queuing networks for EV charging is very limited. It has been studied in \cite{Carvalho2015b}. %We discuss some stability conditions for queuing networks. More specifically, we discuss a numerical stability condition in a queuing framework for EV-charging, and from a more abstract point of view, we go over stability conditions found in the field of \emph{bandwidth-sharing networks}, since our queuing model falls in both categories. 
Here, the authors find by simulation that there is a threshold on the arrival rates of EVs, such that if the actual arrival rate is greater than this threshold, some cars have to wait for increasingly long times to fully charge. The first analytical study is \cite{Christianen2021}, where the authors compare these thresholds on the arrival rates under different power flow models and compute the difference between these rates explicitly as the number of parking lots grows to infinity.

As said earlier, the queuing model used can be viewed as a bandwidth-sharing network. Specifically, EV owners can charge their EVs at parking lots at charging stations. In practice, EVs are served simultaneously, because all the charging stations in the network are connected by underground cables through which power is transmitted from a power source. Thus, all EVs share this power. This feature of sharing limited service capacity with concurrent users can be captured adequately in queuing theory by a resource-sharing network, specifically a bandwidth-sharing network. 

Bandwidth-sharing networks have been successfully applied in many areas, for example in computer systems and communication networks. There is a vast literature on stability conditions in bandwidth-sharing networks; see \cite{Wang2022} for a complete literature review. The question of stability is examined under Markovian assumptions in \cite{Veciana1999,DeVeciana2001}, where the stability is covered for weighted proportionally fair policies. In \cite{Bonald2001a}, these results are generalized to a more general family of policies called weighted $\alpha$-fair policies \cite{Bonald2006a}. Stability results for more general arrival and service distributions have also been derived \cite{Walrand2000,Gromoll2008,Chiang2006}. In particular, \cite{Massoulie2007} has shown the stability of the proportionally fair policy under phase-type service distributions. In a continuous-time Markovian setting, \cite{Shneer2018,Shneer2019c} present stability conditions for a more general class of utility-maximizing allocations, which include weighted $\alpha$-fair allocations. 
  
\paragraph{Comparison of the Linearized Distflow and Distflow model.}
The practical use of the Linearized Distflow compared to the Distflow model is based on the assumption that power losses on cables are typically small. It has been shown experimentally that this only introduces a small relative error, on the order of $1\%$ \cite{Farivar2013}. However, these small relative errors may be exaggerated when used in a complex stochastic process and this is what we examine in this paper. Multiple other numerical studies have been conducted to verify the accuracy and effectiveness of the Linearized Distflow model \cite{Baran1989,Wang2014,Chen2016,Tan2013,Yuan2016,Wang_inverter_2014,Yeh2012a,Li2019,Cao2020}.  

\section{\uppercase{Model description and formulation}}\label{sec:problem_description}
This section describes the main components of the EV-charging model, i.e.; we describe the characteristics of the queuing, the distribution network, and the power flow models. 

\subsection{Queuing model of EV-charging}\label{subsec:queuing_model}
We use a queuing model to study the process of charging EVs in a distribution network. EVs referred to as \emph{jobs} require service. This service is delivered by charging stations, referred to as \emph{servers} and the service being delivered is the power supplied to EVs. At all parking lots, there is a charging station with $K>0$ parking spaces, and each parking space has its own EV charger. %Since we are interested in the maximal feasible arrival rate, we assume that we can simultaneously charge an infinite number of EVs so that we are not limited by space, i.e. a finite number of charging stations or limited numbers of parking lots. At the same time, this is not of the physical realm if we consider, for example, wireless charging.

Thus, in the queuing system, we consider $N$ single-server queues, each having its own arrival stream of jobs. Denote by $\mathbf{X}(t) = (X_1(t),\ldots,X_N(t))$ the vector giving the number of jobs at each queue at time $t$. At all parking lots, all EVs arrive independently according to Poisson processes with rate $\lambda_i, i=1,\ldots,N$ and have independent service requirements which are $Exp(1)$ random variables.  %Second, for the Distflow model, the computat,ion of an arrival rate such that the maximal voltage drop is attained uses an approximation of the voltages under the Distflow model with the same arrival rate at every node.
 %Assumption \ref{assumption:arrival_rate} %is probably not the most pragmatic and realistic assumption, however we make this assumption to 
%keeps the analysis in Section \ref{subsec:distflow} tractable.
%We make this assumption, since our goal is to compare different power flow models. Therefore, we exclude other sources of variability from our model, such as different arrival and demand patterns.

At each queue, all jobs are served simultaneously and start service immediately (there is no queuing). Furthermore, each job receives an equal fraction of the service capacity allocated to a queue. Denote by $\mathbf{\tilde{p}}(t) = (\tilde{p}_1(t),\ldots,\tilde{p}_N(t))$ the vector of service capacities allocated to each queue at time $t$. From now on, for simplicity, we drop the dependence on time $t$ from the notation. For example, we write $X_j$ and $\tilde{p}_j$ instead of $X_j(t)$ and $\tilde{p}_j(t)$.

Service capacities are state-dependent and subject to changes to the current vector $\mathbf{X}=(X_1,\ldots, X_N)$ of number of jobs. For each state of the system, i.e. a given number of EVs charging at each parking lot, we assume that the charging rates $\mathbf{\tilde{p}}$ are the unique solution of the optimization problem
\begin{align}
    \max_{\mathbf{\tilde{p}}} \sum_j X_j\log\bigg(\frac{\tilde{p}_j}{X_j}\bigg),\label{eq:opt_p}
\end{align} which are called proportional fair allocations. For the optimization problem, the feasible region can take many forms and depends heavily on the power flow model that is used. In Section \ref{subsec:power_flow}, we discuss the feasible regions for both power flow models in more detail.

We can then represent the number of electric vehicles charging at every station as an $N$-dimensional continuous-time Markov process. The evolution of the queue at node $j$ is given by
\begin{align*}
X_j(t) \to X_j(t)+1~\text{at rate}~\lambda_j,
\end{align*} and
\begin{align*}
X_j(t) \to X_j(t)-1~\text{at rate}~\tilde{p}_j.
\end{align*}
  
\subsection{Distribution network model}\label{subsec:distribution_network}
The distribution network is modeled as a directed graph $\mathcal{G}=(\mathcal{I},\mathcal{E})$, where we denote by $\mathcal{I} = \{0,\ldots,N\}$ the set of nodes and by $\mathcal{E}$ its set of directed edges, assuming that node $0$ is the root node. We assume that $\mathcal{G}$ has a line topology. Each edge $\epsilon_{j-1,j}\in\mathcal{E}$ represents a line connecting nodes $j-1$ and $j$ where node $j$ is further away from the root node than node $j-1$. Each edge $\epsilon_{j-1,j}\in\mathcal{E}$ is characterized by the impedance $z=r+\mathrm{i} x$, where $r,x\geq 0$ denote the resistance and reactance along the lines, respectively. We make the following natural assumption, given that $r>>x$ in distribution networks \cite{Khatod2006,Tonso2005}.
\begin{assumption}\label{assumption:reactance}
All edges have the same resistance value $r>0$ and reactance value $x=0$.% and reactance value $x$.
\end{assumption} %Alternatively, the corresponding mutual admittance is $y = g+\mathrm{i}b=1/z$. 
 
%We assume that the phase angle between voltages $
%\tilde{V}_i$ and $\tilde{V}_j$ is small in distribution networks \cite{Carvalho2015b}, and hence the phases of $
%\tilde{V}_i$ and $
%\tilde{V}_j$ are approximately the same and can be chosen so that the phasors have zero imaginary components. % In general, there is also the notion of a complex form for the voltage.  However, we assume that voltages have zero imaginary components \cite{Carvalho2015b}. For an explanation of why this is a reasonable assumption, we refer to \cite{Kerstinga}. 
%For $j\in\mathcal{I}$, $\tilde{V}_j$ denotes the real voltage and it emerges that the impedance is zero, thus $z=r$ and all edges have the same resistance value $r$.
%\begin{assumption}

\begin{figure*}[h!]
\begin{center}
\begin{tikzpicture}[main_node/.style={circle,fill=blue!10,minimum size=1em,inner sep=4pt},
feeder_node/.style={ellipse,fill=yellow!10,minimum size=3em,inner sep=4pt},]

    \node[feeder_node][label={$\tilde{V}_0$}] (1) at (0,0) {Root node};
    \node[main_node][label={$\tilde{V}_1$}] (2) at (4,0)  {1};
    \node[main_node][label={$\tilde{V}_i$}] (3) at (8,0) {i};
    \node[main_node][label={$\tilde{V}_N$}] (4) at (12,0) {N};
    \node[below=1.2cm of 2](D){};
    \node[below=1.2cm of 3](E){};
    \node[below=1.2cm of 4](F){};

    \draw (1.4,0) -- (3.7,0); 
    \draw[dashed] (4.3,0) -- (7.7,0)node[midway,above]{$(\tilde{S}_{ij},I_{ij},r)$};
    \draw[dashed] (8.3,0) -- (11.7,0);
    \Edge[Direct,label=$\lambda_1$](D)(2);
    \Edge[Direct,label=$\lambda_i$](E)(3);
    \Edge[Direct,label=$\lambda_N$](F)(4);
    	%(1) .. node[below] {$R$} (2)
    	%(2) edge node[below] {$R$} (3)
    	%(3) edge node[below] {$R$} (4);
\end{tikzpicture}
\end{center}
\caption{Line network with $N$ charging stations and arriving vehicles at rate $\lambda_i$, $i\in\{1,\ldots,N\}$.}
\label{fig:model}
\end{figure*}
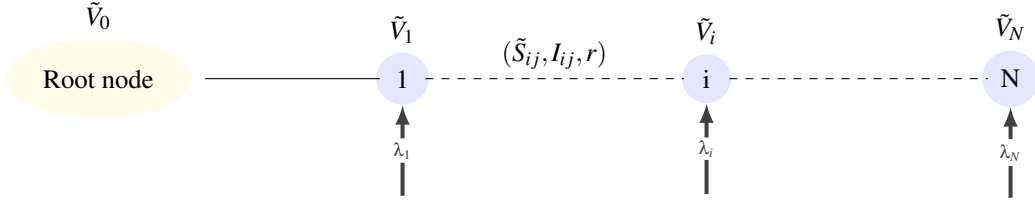

%\end{assumption}
Furthermore, let $\tilde{s}_j = \tilde{p}_j + \mathrm{i} \tilde{q}_j$ be the complex power consumption at node $j$. Here, $\tilde{p}_j$ and $\tilde{q}_j$ denote the active and reactive power consumption at node $j$, respectively. By convention, a positive active (reactive) power term corresponds to consuming active (reactive) power. Since EVs can only consume active power \cite{Carvalho2015b}, it is natural to make the following assumption.
\begin{assumption}\label{assumption:reactive}
The active power $\tilde{p}_j$ is non-negative and the reactive power $\tilde{q}_j$ is zero at all charging stations $j\in\mathcal{I}$.
\end{assumption} Let $\tilde{V}_j$ denote the voltage at node $j$. Given Assumptions \ref{assumption:reactance} and \ref{assumption:reactive}, the voltages at each node $j$, $\tilde{V}_j$, can be chosen to have zero imaginary components \cite{Carvalho2015b,Aveklouris2019}.
%\begin{assumption}
%The reactive power consumption is equal to zero at every charging station, i.e.,
%\begin{align*}
%\tilde{q}_j = 0\ \text{for all}\ j=1,\ldots,N.
%\end{align*}
%\end{assumption}
For each $\epsilon_{j-1,j}\in \mathcal{E}$, let $I_{j-1,j}$ be the complex current and $\tilde{S}_{j-1,j}=\tilde{P}_{j-1,j}+\mathrm{i}\tilde{Q}_{j-1,j}$ be the complex power flowing from node $j-1$ to $j$. Here, $\tilde{P}_{j-1,j}$ and $\tilde{Q}_{j-1,j}$ denote the active and reactive power flowing from node $j-1$ to $j$. 
The model is illustrated in Figure \ref{fig:model}.

\paragraph{Voltage drop constraint} The distribution network constraints, that is in our case only the voltage drop constraint, represent the feasible region of \eqref{eq:opt_p} and are described by a set $\mathcal{C}$. The set $\mathcal{C}$ is contained in an $N$-dimensional vector space and represents feasible power allocations. In our setting, a power allocation is feasible if the maximal voltage drop; i.e., the relative difference between the base voltage $\tilde{V}_0$ and the minimal voltage in all buses between the root node and any other node is bounded by some value $\Delta\in (0,\frac{1}{2}]$. Thus, the distribution network constraints can be described as
\begin{align}
\mathcal{C}:= \left\{\mathbf{\tilde{p}}:\ \frac{\tilde{V}_0-\min_{1\leq j\leq N} \tilde{V}_j}{\tilde{V}_0}\leq \Delta,\quad 0<\Delta\leq \frac{1}{2} \right\}.\label{eq:voltage_drop_C}
\end{align} %In Section \ref{subsec:distflow}, we show there is a technical reason why we use the explicit interval $(0,\frac{1}{2}]$ for $\Delta$. However, this is not out of the physical realm, since in practice, we want the maximal voltage drop to be no more than a small percentage of the base voltage, e.g. $\Delta = 0.05$ or $\Delta=0.1$. 
In Section \ref{subsec:power_flow}, we give more concrete definitions of the constraint set $\mathcal{C}$ for each power flow model. 
\subsection{Power flow models}\label{subsec:power_flow}
We introduce two commonly used models to represent the power flow that is valid for radial systems; i.e., systems where all charging stations have only one (and the same) source of supply. % The Distflow equations fully represent the power flows for a balanced, single-phase equivalent model of a radial network. %One obtains a set of equations that fully represent the power flows in mesh networks by augmenting the DistFlow equations with constraints that ensure that there exist phase angles which, 1) are consistent with the squared voltage magnitudes, |V_i|^2, and power flows, P_ik and Q_ik, and 2) sum to a multiple of 2\pi radians around every cycle.
They are called the \emph{Distflow} and \emph{Linearized Distflow model} \cite{Low2014a,Baran1989}. Both models are valid when the underlying network topology is a tree, which is indeed the case in this paper (as we consider a line topology). 

Given the impedance $r$, %$z$, 
the voltage at the root node $\tilde{V}_0$ and the power consumptions $\tilde{p}_j, j=1,\ldots,N$, both power flow models satisfy three relations. First, we have power balance at each node: for all $j\in\mathcal{I}\backslash\{0\}$,
\begin{align}
\tilde{S}_{j-1,j}-r\left| I_{j-1,j}\right|^2 = \tilde{s}_j + \tilde{S}_{j,j+1}.\label{eq:power_flow_equations}
\end{align} Here, on the one hand, the quantity $r|I_{j-1,j}|^2$ represents line loss so that $\tilde{S}_{j-1,j}-r|I_{j-1,j}|^2$ is the receiving-end complex power at node $j$ from node $j-1$. On the other hand,  the delivering-end complex power is the sum of the consumed power at node $j$ and the complex power flowing from node $j$ to node $j+1$. Second, by Ohm's law, we have for each edge $\epsilon_{j-1,j} \in \mathcal{E}$,
\begin{align}
\tilde{V}_{j-1} - \tilde{V}_j = rI_{j-1,j}\label{eq:ohm}
\end{align} and third, due to the definition of complex power, we have for each edge $\epsilon_{j-1,j} \in \mathcal{E}$,
\begin{align}
\tilde{S}_{j-1,j} = \tilde{V}_{j-1}\overline{I}_{j-1,j}.\label{eq:conjugate}
\end{align} 

However, for a line topology, the distribution network constraints in \eqref{eq:voltage_drop_C}, using \eqref{eq:power_flow_equations}--\eqref{eq:conjugate}, can be rewritten for both power flow models \cite{Christianen2021}. Following the approach in \cite{Christianen2021}, the distribution network constraints under the Distflow model $\mathcal{C}^D$ reduce to
\begin{align}
\mathcal{C}^D:= \left\{\mathbf{\tilde{p}}:\ \tilde{V}_0^D \leq \frac{1}{1-\Delta}, \right\}, \quad 0<\Delta\leq \frac{1}{2},\label{eq:voltage_drop_C_D}
\end{align} where $\tilde{V}_0^D$ can be found recursively by
\begin{align}
    \tilde{V}_{N-1}^D & = 1+r\tilde{p}_N,\\
    \tilde{V}_{j-1}^D & = 2\tilde{V}_j^D-\tilde{V}_{j+1}^D+\frac{r\tilde{p}_j}{\tilde{V}_j^D},\quad j=1,\ldots,N-1.
\end{align} In the Linearized Distflow model, it is assumed that the active power losses $r|I_{j-1,j}|^2$ are much smaller than the active power flows $\tilde{P}_{j-1,j}$. In other words, the Linearized Distflow model neglects the loss terms associated with the squared current magnitudes $|I_{j-1,j}|^2$. In that case, the distribution network constraints under the Linearized Distflow $\mathcal{C}^{LD}$ reduce to
\begin{align}
    \mathcal{C}^{LD} :=  \left\{\mathbf{\tilde{p}}:\ 2r\sum_{j=1}^{N}\sum_{k=j}^{N}\tilde{p}_{k} \leq \frac{\Delta(2-\Delta)}{(1-\Delta)^2}, \right\}, \quad 0<\Delta\leq \frac{1}{2}.\label{eq:voltage_drop_C_LD}
\end{align}

\section{\uppercase{Numerical study}}\label{sec:numerical_study}
In the previous section, we discussed our model for the EV-charging process. Here, we obtain general insights into the performance of the model by simulation on a large range of parameter settings. Moreover, this allows us to compare the behavior of the model under the Distflow model and the Linearized Distflow model. We vary the total arrival rate to the network and the distribution of the total arrival rate to different parking lots. We focus on the effect of the mean number of EVs in the network and the mean charging time of an EV, possibly at each parking lot.

\subsection{Critical arrival rate} 
\begin{figure*}
\centering
\begin{subfigure}{.5\textwidth}
  \centering
  \includegraphics[width=\linewidth]{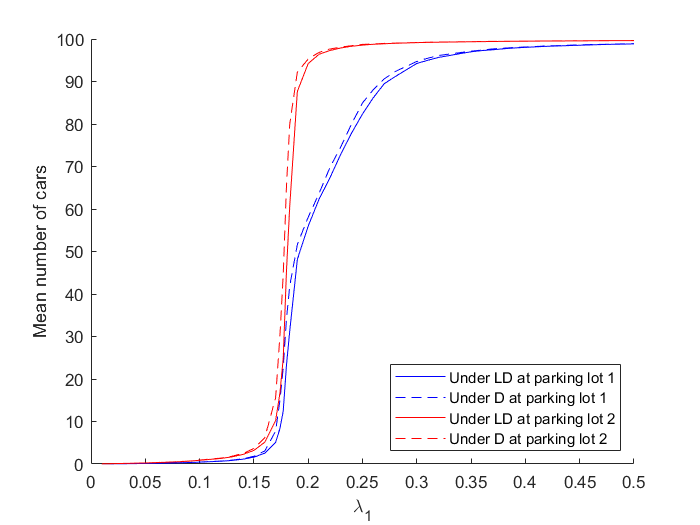}
  \caption{Mean number of cars}
  \label{fig:sub1}
\end{subfigure}%
\begin{subfigure}{.5\textwidth}
  \centering
  \includegraphics[width=\linewidth]{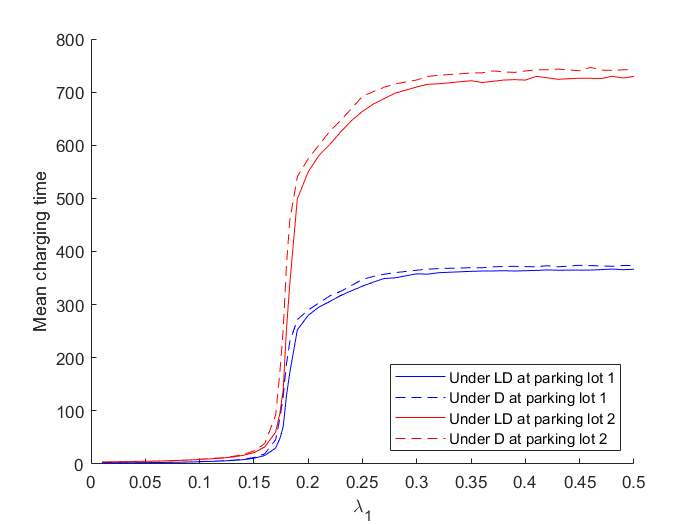}
  \caption{Mean charging time}
  \label{fig:sub2}
\end{subfigure}
\caption{Performance measures vs. the individual arrival rate per parking lot for the Distflow model (dashed) and the Linearized Distflow model (solid) at parking lot 1 (blue) and at parking lot 2 (red).}
\label{fig:critical_arrival_rate}
\end{figure*}

To control the network, we observe that there is a \emph{critical} arrival rate $\lambda_c$ ($\approx 0.18$) if the arrival rate to each parking lot is assumed to be the same. At every parking lot, the mean number of EVs and the mean charging time of an EV seem to explode as soon as the actual arrival rate is greater than the critical arrival rate. See Figure \ref{fig:critical_arrival_rate}, where the mean number of EVs and the mean charging time of an EV are plotted versus the individual arrival rate to each parking lot. We observe the mean number of EVs and the mean charging time of an EV for both power flow models at each parking lot. We fix the number of parking lots at $N=2$, the resistance for each cable at $r=0.1$, the maximal capacity at each parking lot at $K=100$ and the parameter to control the voltage drop at $\Delta = 0.05$ for the Distflow model (dashed) and the Linearized Distflow model (solid) at parking lot 1 (blue) and at parking lot 2 (red). Although the total arrival rate is varied, the solid curves are close to the dashed lines for the two parking lots. Moreover, for all arrival rates, the solid curves are below the dashed curves. This is to be expected; in \cite{Christianen2021}, the authors observed that the
Linearized Distflow power flow model allows for too optimistic arrival rates since the Linearized Distflow model overestimates the voltages in comparison with the voltages given by the Distflow model. From this, we observe that the allocated power to each parking lot is higher under the Linearized Distflow than the allocated power under the Distflow model. Higher allocated power means faster charging. Hence, EVs leave the parking lots faster and the mean number of EVs charging is lower.

\begin{figure}
    \centering
    \includegraphics[width=\linewidth]{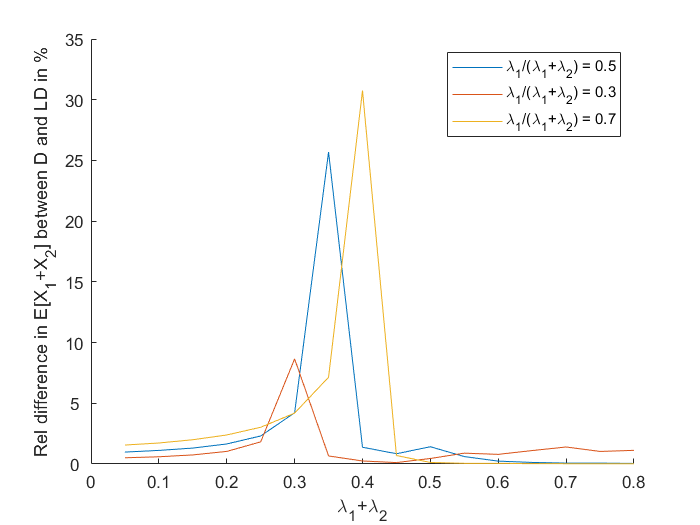}
    \caption{Relative difference in mean number of EVs between both power flow models}
    \label{fig:rel_difference_D_and_LD}
\end{figure}

This statement is reinforced by Figure \ref{fig:rel_difference_D_and_LD}, where the relative difference in the mean number of EVs in the network in percentages between the Distflow model and the Linearized Distflow model is plotted versus the total arrival rate to the parking lots. It is still assumed that $N=2$, $r=0.1$, $K=100$ and $\Delta=0.05$ for all of these lines. For the blue curve, we have equal arrival rates for all parking lots. From Figure \ref{fig:rel_difference_D_and_LD}, it is apparent that the relative difference in the mean number of EVs between both power flow models is below $5\%$ for almost all total arrival rates to the network. However, when the total arrival rate is close to two times the critical arrival rate $\lambda_c (\approx 0.36)$, we have already seen in Figure \ref{fig:critical_arrival_rate} that the mean number of cars for both power flow models seem to explode and that this happens slightly earlier for the Distflow model than for the Linearized Distflow model. This causes the high relative difference in mean number of EVs for both power flow models.

\subsection{Variability of the distribution of the total arrival rate}

The previous section brings us naturally to the discussion of adding variability to the distribution of the total arrival rate. Instead of assuming equal arrival rates for all parking lots, we vary the fraction of EVs that arrive at each parking lot for a wide range of the total arrival rate to the network.

\begin{figure}
    \centering  
      \includegraphics[width=\linewidth]{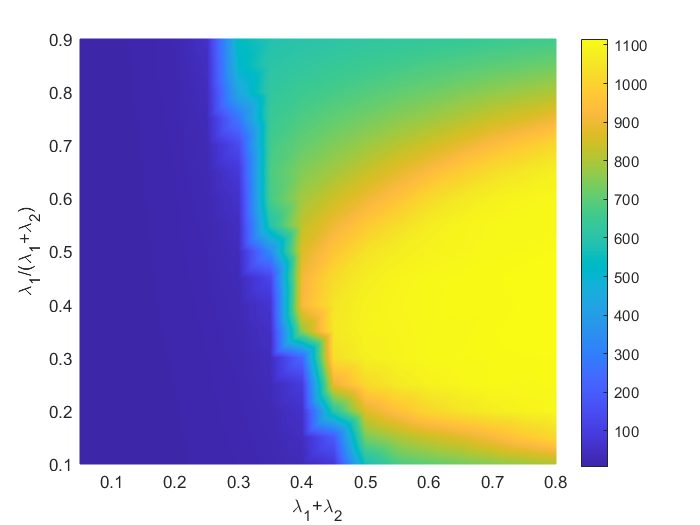}
  \caption{Mean number of cars}
  \label{fig:heatmap}
\end{figure}

\begin{figure*}
\centering
\begin{subfigure}{.5\textwidth}
  \centering
\includegraphics[width=\linewidth]{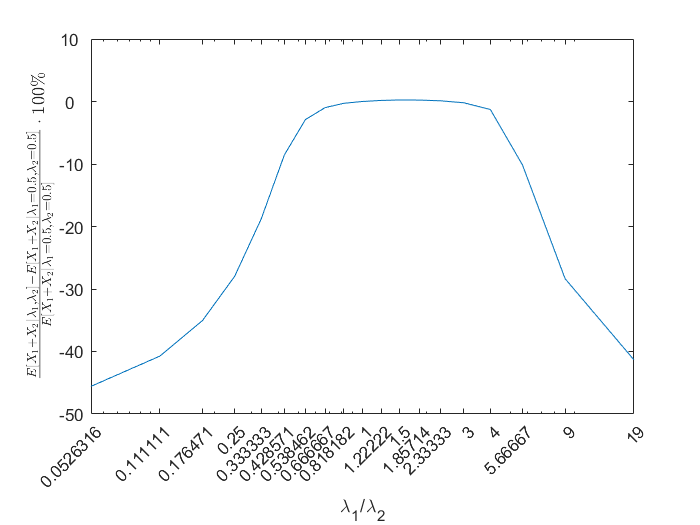}
  \caption{Relative difference for $\lambda_1+\lambda_2 = 0.8$}
  \label{fig:rel_difference}
\end{subfigure}%
\begin{subfigure}{.5\textwidth}
  \centering

  \includegraphics[width=\linewidth]{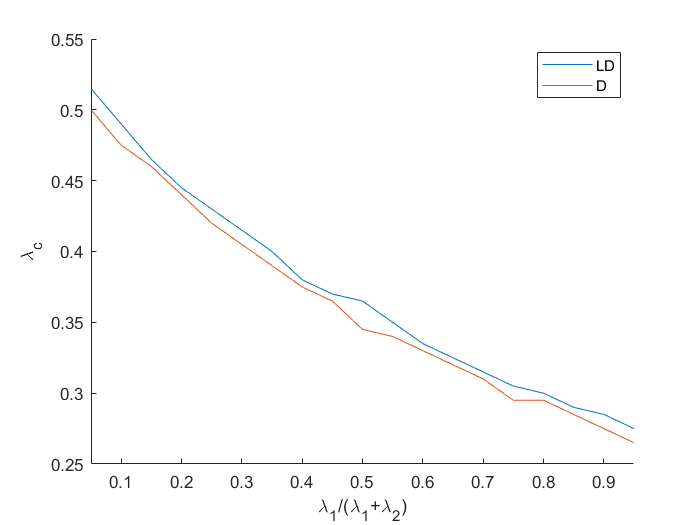}
    \caption{Approximation of critical arrival rates}
    \label{fig:approx_crit_arrival_rtates}
\end{subfigure}
\caption{}
\label{fig:heatmap_rel_difference}
\end{figure*}

For all combinations of total arrival rates and fractions of EVs that arrive at each parking lot, the heat map of the mean number of cars has an interesting structure. In Figure \ref{fig:heatmap}, we show the total mean number of EVs in the network as a function of the total arrival rate and the fractions of EVs that arrive at each parking lot. As we see in Figure \ref{fig:heatmap}, the mean number of EVs has a non-symmetric structure. When we increase the fraction of EVs that arrive at parking lot 1 (and thus decrease the fraction of EVs that arrive at parking lot 2) from the situation of equal arrival rates, the total mean number of EVs in the network decreases faster than the total number of EVs in the network decreases when we increase the fraction of EVs that arrive at parking lot 2. This is natural given the total available power that can be allocated to each parking lot. Due to the power loss on the cables, the available power that can be allocated to parking lot 1 is approximate twice the available power that can be allocated to parking lot 2. We compare the following two situations; one where we have a given number of EVs charging at parking lot 1 and no EVs charging at parking lot 2 (which corresponds to a situation where the fraction of EVs that arrive at parking lot 1 is high), and one where we have the same given number of EVs charging at parking lot 2 and no EVs charging at parking lot 1 (which corresponds to a situation where the fraction of EVs that arrive to parking lot 2 is high). Since the allocated power to parking lot 1 in the first situation is higher than the allocated power to parking lot 2 in the second situation, the mean number of EVs at parking lot 1 tends to be smaller than the number of EVs at parking lot 2.   Moreover, if we consider a fixed total arrival rate, for example, $\lambda_1+\lambda_2 = 0.8$, the variability of the distribution of the total arrival rate has a small influence on the mean number of EVs in the network for a wide range of the ratio of arrival rates of both parking lots. In Figure \ref{fig:rel_difference}, where the relative difference between the total mean number of EVs given any distribution of the total arrival over both parking lots and the total mean number of EVs given equal arrival rates to both parking lots is plotted, we observe that this relative difference is below $5\%$ if the fraction of EVs that arrive to parking lot 1 range between $20\%$ and $60\%$.

% Two queues explode (yellow) vs one queue explodes (green)
Another observation on the heat map of the mean number of cars is that there is a clear distinction between networks that have or have not reached their capacity. In the blue region, the mean number of EVs is relatively low. However, in the green and yellow regions, the mean number of EVs is relatively high and close to its maximal capacity. Moreover, the green region indicates a network where the number of EVs charging at either parking lot 1 or 2 explodes, or in other words, that either parking lot 1 or 2 reaches its capacity. The yellow region indicates a network where the number of EVs at both parking lots explodes.

In Figure \ref{fig:approx_crit_arrival_rtates}, we also plot a rough approximation of the evolution of the critical arrival rates under both power flow models as we vary the fraction of EVs that arrive at parking lot 1 that we obtained as follows. As can be seen in Figure \ref{fig:critical_arrival_rate}, there is a steep increase in the mean number of EVs and the mean charging time of an EV at a certain point. Using a fine grid search, we find the arrival rate for which the absolute difference between subsequent measurements in the mean number of cars is the largest. In Figure \ref{fig:approx_crit_arrival_rtates}, we observe that as the fraction of EVs that arrive at parking lot 1 increases, the critical arrival rates for both power flow models decrease. Furthermore, the critical arrival rates for the Distflow model are always below the critical arrival rates for the Linearized Distflow model; a behavior that we observed in the setting of equal arrival rates for both parking lots already.

Adding variability to the distribution of the total arrival also influences the relative difference of the total mean number of EVs in the network between both power flow models. This effect is observed in Figure \ref{fig:rel_difference_D_and_LD}, where we see that for a larger fraction of EVs that arrive at parking lot 1, the maximal relative difference between the two power flow models increases. Furthermore, also for other distributions of the total arrival rate than equal arrival rates to both parking lots (red and yellow curves), we observe that the relative difference is below $5\%$ for almost all total arrival rates to the network, except around those values for the total arrival rate where we turn from a network with a relatively low mean number of EVs to a network with a relatively high number of EVs. 

In summary, the performance of both power flow models is approximately the same. Simulation results show that for a wide range of total arrival rates to the network the relative difference in the mean number of EVs between the Distflow model and the Linearized Distflow model is below $5\%$. Furthermore, the critical arrival rates under both power flow models are close to each other. Moreover, we can say that the variability in the distribution of the total arrival to parking lots, as long as it is not too large, does not influence the performance of the network much, in the sense that the mean number of EVs and the mean charging time of an EV are comparable to the mean number of EVs and the mean charging time of an EV in the case where the arrival rates to all parking lots are the same, respectively.

\typeout{}
\bibliographystyle{apalike}
{\small
\bibliography{library}}

\begin{thebibliography}{}

\bibitem[Aveklouris et~al., 2019]{Aveklouris2019}
Aveklouris, A., Vlasiou, M., and Zwart, B. (2019).
\newblock {A stochastic resource-sharing network for electric vehicle
  charging}.
\newblock {\em IEEE Transactions on Control of Network Systems},
  6(3):1050--1061.

\bibitem[Baran and Wu, 1989a]{Baran1989a}
Baran, M.~E. and Wu, F.~F. (1989a).
\newblock {Optimal capacitor placement on radial distribution systems}.
\newblock {\em IEEE Transactions on Power Delivery}, 4(1):725--734.

\bibitem[Baran and Wu, 1989b]{Baran1989}
Baran, M.~E. and Wu, F.~F. (1989b).
\newblock {Optimal Sizing of Capacitors Placed On a Radial Distribution
  System}.
\newblock {\em IEEE Transactions on Control of Network Systems}, 4(1):735--743.

\bibitem[Bonald and Massouli{\'{e}}, 2001]{Bonald2001a}
Bonald, T. and Massouli{\'{e}}, L. (2001).
\newblock {Impact of fairness on internet performance}.
\newblock {\em Performance Evaluation Review}, 29(1):82--91.

\bibitem[Bonald and Prouti{\`{e}}re, 2006]{Bonald2006a}
Bonald, T. and Prouti{\`{e}}re, A. (2006).
\newblock {Flow-level stability of utility-based allocations for non-convex
  rate regions}.
\newblock {\em 2006 IEEE Conference on Information Sciences and Systems, CISS
  2006 -- Proceedings}, pages 327--332.

\bibitem[Cao et~al., 2019]{Cao2020}
Cao, Y., Wei, W., Wang, J., Mei, S., Shafie-Khah, M., and Catalao, J.~P.
  (2019).
\newblock {Capacity Planning of Energy Hub in Multi-carrier Energy Networks: A
  Data-driven Robust Stochastic Programming Approach}.
\newblock {\em IEEE Power and Energy Society General Meeting},
  2019-Augus(1):3--14.

\bibitem[Carvalho et~al., 2015]{Carvalho2015b}
Carvalho, R., Buzna, L., Gibbens, R., and Kelly, F. (2015).
\newblock {Critical behaviour in charging of electric vehicles}.
\newblock {\em New Journal of Physics}, 17(9):95001.

\bibitem[Chen et~al., 2016]{Chen2016}
Chen, C., Wang, J., Qiu, F., and Zhao, D. (2016).
\newblock {Resilient Distribution System by Microgrids Formation after Natural
  Disasters}.
\newblock {\em IEEE Transactions on Smart Grid}, 7(2):958--966.

\bibitem[Chiang et~al., 2006]{Chiang2006}
Chiang, M., Shah, D., and Tang, A. (2006).
\newblock {Stochastic Stability Under Network Utility Maximization : General
  File Size Distribution}.
\newblock {\em 44st Annual Allerton Conference on Communication, Control, and
  Computing}, pages 1--22.

\bibitem[Christianen et~al., 2022]{Christianen2021}
Christianen, M., Cruise, J., Janssen, A., Shneer, S., Vlasiou, M., and Zwart,
  B. (2022).
\newblock {Comparison of stability regions for a line distribution network with
  stochastic load demands}.
\newblock {\em Preprint available at: https://arxiv.org/abs/2201.06405}.

\bibitem[de~Hoog et~al., 2014]{DeHoog2014}
de~Hoog, J., Muenzel, V., Jayasuriya, D.~C., Alpcan, T., Brazil, M., Thomas,
  D.~A., Mareels, I., Dahlenburg, G., and Jegatheesan, R. (2014).
\newblock {The importance of spatial distribution when analysing the impact of
  electric vehicles on voltage stability in distribution networks}.
\newblock {\em Energy Systems}, 6(1):63--84.

\bibitem[{De Veciana} et~al., 2001]{DeVeciana2001}
{De Veciana}, G., Lee, T.~J., and Konstantopoulos, T. (2001).
\newblock {Stability and performance analysis of networks supporting elastic
  services}.
\newblock {\em IEEE/ACM Transactions on Networking}, 9(1):2--14.

\bibitem[Deb et~al., 2018]{Deb2018}
Deb, S., Tammi, K., Kalita, K., and Mahanta, P. (2018).
\newblock {Impact of electric vehicle charging station load on distribution
  network}.
\newblock {\em Energies}, 11(1):1--25.

\bibitem[Farivar et~al., 2013]{Farivar2013}
Farivar, M., Chen, L., and Low, S. (2013).
\newblock {Equilibrium and dynamics of local voltage control in distribution
  systems}.
\newblock {\em Proceedings of the IEEE Conference on Decision and Control},
  pages 4329--4334.

\bibitem[Gromoll and Williams, 2008]{Gromoll2008}
Gromoll, H.~C. and Williams, R.~J. (2008).
\newblock {Fluid Model for a Data Network with $\alpha$ -Fair Bandwidth Sharing
  and General Document Size Distributions: Two Examples of Stability}.
\newblock 4:253--265.

\bibitem[Hoogsteen et~al., 2017]{Hoogsteen2017}
Hoogsteen, G., Molderink, A., Hurink, J.~L., Smit, G.~J., Kootstra, B., and
  Schuring, F. (2017).
\newblock {Charging electric vehicles, baking pizzas, and melting a fuse in
  Lochem}.
\newblock {\em CIRED - Open Access Proceedings Journal}, 2017(1):1629--1633.

\bibitem[Huang et~al., 2013]{Huang2013}
Huang, H., Chung, C.~Y., Chan, K.~W., and Chen, H. (2013).
\newblock {Quasi-Monte Carlo based probabilistic small signal stability
  analysis for power systems with plug-in electric vehicle and wind power
  integration}.
\newblock {\em IEEE Transactions on Power Systems}, 28(3):3335--3343.

\bibitem[Kersting, 2018]{Kerstinga}
Kersting, W. (2018).
\newblock {\em {Distribution System Modeling and Analysis}}.
\newblock CRC Press, fourth edition.

\bibitem[Khatod et~al., 2006]{Khatod2006}
Khatod, D.~K., Pant, V., and Sharma, J. (2006).
\newblock {A novel approach for sensitivity calculations in the radial
  distribution system}.
\newblock {\em IEEE Transactions on Power Delivery}, 21(4):2048--2057.

\bibitem[Li et~al., 2019]{Li2019}
Li, R., Wei, W., Mei, S., Hu, Q., and Wu, Q. (2019).
\newblock {Participation of an Energy Hub in Electricity and Heat Distribution
  Markets: An MPEC Approach}.
\newblock {\em IEEE Transactions on Smart Grid}, 10(4):3641--3653.

\bibitem[Low, 2014]{Low2014a}
Low, S.~H. (2014).
\newblock {Convex relaxation of optimal power flow - Part {\{}I{\}}:
  Formulations and equivalence}.
\newblock {\em IEEE Transactions on Control of Network Systems}, 1(1):15--27.

\bibitem[Massouli{\'{e}}, 2007]{Massoulie2007}
Massouli{\'{e}}, L. (2007).
\newblock {\em {Structural properties of proportional fairness: Stability and
  insensitivity}}, volume~17.

\bibitem[Mo and Walrand, 2000]{Walrand2000}
Mo, J. and Walrand, J. (2000).
\newblock {Fair End-to-End Window-Based Congestion Control}.
\newblock {\em IEEE/ACM Transactions on Networking}, 8(5):556--567.

\bibitem[Molzahn and Hiskens, 2019]{Molzahn2019}
Molzahn, D. and Hiskens, I. (2019).
\newblock {A Survey of Relaxations and Approximations of the Power Flow
  Equations}.
\newblock {\em A Survey of Relaxations and Approximations of the Power Flow
  Equations}, 4(1):1--221.

\bibitem[Shneer and Stolyar, 2018]{Shneer2018}
Shneer, S. and Stolyar, A. (2018).
\newblock {Stability and moment bounds under utility-maximising service
  allocations: finite and infinite networks}.
\newblock {\em arXiv}, pages 1--26.

\bibitem[Shneer and Stolyar, 2019]{Shneer2019c}
Shneer, S. and Stolyar, A. (2019).
\newblock {Stability conditions for a decentralised medium access algorithm:
  single- and multi-hop networks}.
\newblock {\em Queueing Systems}, 94(1):109--128.

\bibitem[Tan et~al., 2013]{Tan2013}
Tan, S., Xu, J.~X., and Panda, S.~K. (2013).
\newblock {Optimization of distribution network incorporating distributed
  generators: An integrated approach}.
\newblock {\em IEEE Transactions on Power Systems}, 28(3):2421--2432.

\bibitem[Tonso et~al., 2005]{Tonso2005}
Tonso, M., Morren, J., {De Haan}, S.~W., and Ferreira, J.~A. (2005).
\newblock {Variable inductor for voltage control in distribution networks}.
\newblock {\em 2005 European Conference on Power Electronics and Applications},
  2005.

\bibitem[Ul-Haq et~al., 2015]{Ul-Haq2015}
Ul-Haq, A., Cecati, C., Strunz, K., and Abbasi, E. (2015).
\newblock {Impact of electric hehicle charging on voltage unbalance in an urban
  distribution network}.
\newblock {\em Intelligent Industrial Systems}, 1(1):51--60.

\bibitem[van Westering and Hellendoorn, 2020]{VanWestering2020}
van Westering, W. and Hellendoorn, H. (2020).
\newblock {Low voltage power grid congestion reduction using a community
  battery: Design principles, control and experimental validation}.
\newblock {\em International Journal of Electrical Power and Energy Systems},
  114:105349.

\bibitem[Veciana et~al., 1999]{Veciana1999}
Veciana, D.~E., Lee, T., and Konstantopoulos, T. (1999).
\newblock {Stability and performance analysis of network supporting services
  with rate control -- could the Internet be unstable?}
\newblock {\em Proceedings of IEEE Infocom}.

\bibitem[Wang et~al., 2022]{Wang2022}
Wang, W., Maguluri, S.~T., Srikant, R., and Ying, L. (2022).
\newblock {Heavy-Traffic Insensitive Bounds for Weighted Proportionally Fair
  Bandwidth Sharing Policies}.
\newblock {\em Mathematics of Operations Research}.

\bibitem[Wang et~al., 2014a]{Wang2014}
Wang, Z., Chen, B., Wang, J., Kim, J., and Begovic, M.~M. (2014a).
\newblock {Robust optimization based optimal DG placement in microgrids}.
\newblock {\em IEEE Transactions on Smart Grid}, 5(5):2173--2182.

\bibitem[Wang et~al., 2014b]{Wang_inverter_2014}
Wang, Z., Chen, H., Wang, J., and Begovic, M. (2014b).
\newblock {Inverter-less hybrid voltage/var control for distribution circuits
  with photovoltaic generators}.
\newblock {\em IEEE Transactions on Smart Grid}, 5(6):2718--2728.

\bibitem[Yeh et~al., 2012]{Yeh2012a}
Yeh, H.~G., Gayme, D.~F., and Low, S.~H. (2012).
\newblock {Adaptive VAR control for distribution circuits with photovoltaic
  generators}.
\newblock {\em IEEE Transactions on Power Systems}, 27(3):1656--1663.

\bibitem[Yuan et~al., 2016]{Yuan2016}
Yuan, W., Wang, J., Qiu, F., Chen, C., Kang, C., and Zeng, B. (2016).
\newblock {Robust Optimization-Based Resilient Distribution Network Planning
  Against Natural Disasters}.
\newblock {\em IEEE Transactions on Smart Grid}, 7(6):2817--2826.

\bibitem[Zhang et~al., 2016]{Zhang2016}
Zhang, Y., Song, X., Gao, F., and Li, J. (2016).
\newblock {Research of voltage stability analysis method in distribution power
  system with plug-in electric vehicle}.
\newblock {\em Asia-Pacific Power and Energy Engineering Conference, APPEEC},
  2016-Decem(51177152):1501--1507.

\end{thebibliography}

\end{document}